\documentclass[12pt]{amsart} 

\usepackage[english]{babel}
\usepackage[T1]{fontenc}

\usepackage{amssymb} 
\usepackage{amsmath} 
\usepackage{fullpage} 
\usepackage{url}

\newtheorem{theo}{Theorem}
\newtheorem{prop}{Proposition}
\newtheorem{lemma}{Lemma}

\theoremstyle{definition}
\newtheorem{defi}{Definition}
\theoremstyle{remark}
\newtheorem{rem}{Remark}

\newtheorem{result}{Result}
\def\Er{{\mathbb E}}

\def\Pr{{\mathbb P}}
\def\Qr{{\mathbb Q}}
\def\Rr{{\mathbb R}}

\def\Ac{{\mathcal{A}}}
\def\Bc{{\mathcal{B}}}
\def\Cc{{\mathcal{C}}}
\def\Dc{{\mathcal{D}}}
\def\Ec{{\mathcal{E}}}
\def\Fc{{\mathcal{F}}}
\def\Gc{{\mathcal{G}}}

\def\Mc{{\mathcal{M}}}
\def\Nc{{\mathcal{N}}}

\def\one{{\rm \bf 1}}
\def\({\left(}     
\def\){\right)}    
\def\[{\left[}     
\def\]{\right]}

\def\as{{\frenchspacing a.s.}~}
\begin{document}
\title{Conditionally atomless extensions of sigma algebras}
\author{Freddy Delbaen}
\address{Departement f\"ur Mathematik, ETH Z\"urich, R\"{a}mistrasse
101, 8092 Z\"{u}rich, Switzerland}
\address{Institut f\"ur Mathematik,
Universit\"at Z\"urich, Winterthurerstrasse 190,
8057 Z\"urich, Switzerland}
\date{First version March 2020, this version \today}

\begin{abstract}
We give two equivalent definitions of sigma algebras that are atomless conditionally to a smaller sigma algebra. \end{abstract}

\maketitle

\section{Notation}

In this paper  \footnote{ This paper is to be seen as an exercise in measure theory.  It will not be sent to a math. journal} we will work with a probability space equipped with three sigma algebras $(\Omega,\Fc_0\subset \Fc_1\subset \Fc_2,\Pr)$.  The sigma algebra $\Fc_0$ is supposed to be trivial $\Fc_0=\{\emptyset,\Omega\}$ whereas the sigma algebra $\Fc_2$ is supposed to express innovations with respect to $\Fc_1$.  Since we do not put topological properties on the set $\Omega$ we will make precise definitions later that do not use conditional probability kernels.  But essentially we could say that we suppose that conditionally on $\Fc_1$ the probability $\Pr$ is atomless on $\Fc_2$. We will show that such an hypothesis implies that there is an atomless sigma algebra $\Bc\subset \Fc_2$ that is independent of $\Fc_1$. In some (under topological hypotheses on $\Omega,\Fc_1,\Fc_2$) cases the conditional expectation with respect to $\Fc_1$ is given by integration with respect to a kernel. We will use the notation $K$ for such a kernel. More precisely: the mapping $K\colon \Omega\times \Fc_2\rightarrow \Rr_+$ satisfies

\begin{enumerate} 
\item For almost every $\omega\in \Omega$, the mapping $K(\omega,.):\Fc_2\rightarrow [0,1]$ is a probability.  It is no restriction to suppose that this property holds for every $\omega\in \Omega$.
\item For each $A\in \Fc_2$, the mapping $K(.,A)\colon \Omega\rightarrow [0,1]$ is $\Fc_1$ measurable.
\item For each $\xi\in L^1(\Omega,\Fc_2,\Pr)$ we have that almost surely
$$
\Er[\xi\mid\Fc_1](\omega)=\int \xi(\tau)\,K(\omega,d\tau).
$$
\end{enumerate}
The existence of such a kernel is not always easy to verify.  Sometimes it is part of the model that is studied. Applying the property above and integrating with respect to $\Pr$  gives
$$
\Pr[A]=\int_{\Omega}
\Pr[d\omega]K(\omega,A).
$$
Or for general $\xi\in L^1$:
$$
\int_{\Omega}\Pr[dx]\xi(x)=\int_\Omega \Pr[d\omega]\int_\Omega \xi(\tau)\,K(\omega,d\tau).
$$
Part of the results were developed and used in my paper on commonotonicity, see \cite{del-com}

\section{Atomless Extension of Sigma Algebras}

\begin{defi}  We say that $\Fc_2$ is atomless conditionally to $\Fc_1$ if the following holds.  For every  $A\in \Fc_2$  there exists a set $B\subset A$, $B\in\Fc_2$, such that $0< \Er[\one_B\mid\Fc_1]<\Er[\one_A\mid\Fc_1]$ on the set $\{\Er[\one_A\mid\Fc_1]>0\}$.
\end{defi}
In  case the conditional expectation could be calculated with a -- under extra topological conditions -- regular probability kernel, say $K(\omega, A)$, then the above definition is a measure theoretic way of saying that the probability measure $K(\omega, .)$ is atomless for almost every $\omega\in\Omega$.  This equivalence will be the topic of the next section.
x\begin{theo} $\Fc_2$ is atomless conditionally to $\Fc_1$ if for every $A\in \Fc_2$, $\Pr[A]>0$, there is $B\subset A$ such that 
$$\Pr\[0< \Er[\one_B\mid\Fc_1]<\Er[\one_A\mid\Fc_1]\]>0.$$
\end{theo}
{\bf Proof} The proof is a standard exhaustion argument.  For completeness we give a proof.   Let $\Dc$ be the collection of $\Fc_1-$measurable sets:
$$
\Dc=\left\{ \{0<\Er[\one_B\mid\Fc_1]<\Er[\one_A\mid\Fc_1]\}\mid B\subset A     \right\}
$$
We show that there is a biggest set in $\Dc$ and this set must then be equal to 
$\{\Er[\one_A\mid\Fc_1]>0\}$.  To show that there is a biggest set in $\Dc$ it is sufficient to show that $\Dc$ is stable for countable unions.  Let $D_n$ be a sequence in $\Dc$ and suppose that for each $n$ we have a set $B_n\subset A$ such that
$D_n=\{ 0<\Er[\one_{B_n}\mid\Fc_1]<\Er[\one_A\mid\Fc_1]\}$.  Now take
$$
B=\cup_n \( B_n\cap \(D_n\setminus \(\cup_{k\le n-1}D_k\)\)  \).
$$
It is easy to check that $\{ 0<\Er[\one_{B}\mid\Fc_1]<\Er[\one_A\mid\Fc_1]\}=\cup_nD_n$ and therefore $\cup_n D_n\in \Dc$.  Let $D=\{ 0<\Er[\one_{B}\mid\Fc_1]<\Er[\one_A\mid\Fc_1]\}$ be a maximum in $\Dc$. Suppose now that
$\Pr[\{\Er[\one_A\mid\Fc_1]>0\}\setminus D  ]>0$.  This implies that $\Pr[A\setminus D]>0$. According to the hypothesis of the theorem, there will be a set $B'\subset (A\setminus D)$ with $D'\subset\{ 0<\Er[\one_{B'}\mid \Fc_1]<\Er[\one_{A\setminus D}\mid \Fc_1]  \} $ and non-negligible.  Since $D\cup D'\in \Dc$ and $D\cap D'=\emptyset$, the element $D$ is not a maximum, a contradiction.

The main result of this section is the following
\begin{theo} $\Fc_2$ is atomless conditionally to $\Fc_1$ if and only if there exists an atomless sigma algebra $\Bc\subset \Fc_2$ that is independent of $\Fc_1$.
\end{theo}
The ``if" part is easy but requires some continuity argument.  Because $\Bc$ is atomless, there is a $\Bc$-measurable, $[0,1]$ uniformly distributed random variable $U$. The sets $B_t=\{U\le t\}, 0\le t \le 1$ form an increasing family of sets with $\Pr[B_t]=t$. Let $A\in \Fc_2$  and let $F=\{ 0 < \Er[\one_A\mid \Fc_1]\}$.  We may suppose that $\Pr[F]>0$ since otherwise there is nothing to prove. We will show that there is $t\in]0,1[$ with $\Pr\[ 0 < \Er[\one_{A\cap B_t} \mid \Fc_1]< \Er[\one_A\mid \Fc_1]    \] > 0$.  According to the previous theorem, $\Fc_2$ is conditionally atomless with respect to $\Fc_1$.  Obviously for $0\le s\le t \le 1$ we have, by independence of $\Bc$ and $\Fc_1$:
$$
\Vert   \Er[\one_{A\cap B_t} \mid \Fc_1] -  \Er[\one_{A\cap B_s} \mid \Fc_1]\Vert_\infty \le \Vert   \Er[\one_{B_t\setminus B_s} \mid \Fc_1]\Vert_\infty = t-s.
$$
It follows that there is a set of measure $1$, say $\Omega'$,  such that for all $s\le t$, rational, $$| \Er[\one_{A\cap B_t} \mid \Fc_1] -  \Er[\one_{A\cap B_s} \mid \Fc_1] | \le t-s$$ on $\Omega'$.  For $\omega\in \Omega'$ we can extend the function $$\{q\in [0,1] \mid q \text{ rational }\}\rightarrow \Er[\one_{A\cap B_q}\mid \Fc_1](\omega)$$ to a continuous function on $[0,1]$.  The resulting continuous extension then represents $\(\Er
[\one_{A\cap B_t} \mid \Fc_1]\)_{t}$. For $t=0$ we have zero and for $t=1$ we find $\Er[\one_A\mid \Fc_1]$. Because for $\omega\in\Omega'$, the trajectories are continuous,  a simple application of Fubini's theorem shows  that the real valued function
$$
t\rightarrow \Pr\[ 0 < \Er[\one_{A\cap B_t} \mid \Fc_1]< \Er[\one_A\mid \Fc_1]    \] 
$$ 
becomes strictly positive for some $t$.  For completeness let us now give the details of the application of Fubini's theorem.  Suppose on the contrary that for all $t\in [0,1]$ we have 
$$
\Pr\[ 0 < \Er[\one_{A\cap B_t} \mid \Fc_1]< \Er[\one_A\mid \Fc_1]    \] =0.
$$
Then on the product space $[0,1]\times \Omega'$ we find that the (clearly measurable) set $$\{(t,\omega)\mid 0 < \Er[\one_{A\cap B_t} \mid \Fc_1](\omega)< \Er[\one_A\mid \Fc_1] (\omega)\}$$
has $m\times\Pr$ measure zero ($m$ denotes Lebesgue measure).  By Fubini's theorem we have that for almost all $\omega\in\Omega'$, the set
$$
\{t \mid 0 < \Er[\one_{A\cap B_t} \mid \Fc_1](\omega)< \Er[\one_A\mid \Fc_1] (\omega)\}
$$
must have Lebesgue measure zero.  However, for $\omega\in \Omega'$, this contradicts the continuity of the mapping
$$
t\rightarrow \Er[\one_{A\cap B_t} \mid \Fc_1](\omega).
$$
The proof of the ``only if" part is broken down in several steps.  We will without further notice, always suppose that $\Fc_2$ is atomless conditionally to $\Fc_1$.
\begin{lemma} Suppose $A\in\Fc_1$ and $C\subset A$ is such that $\Er[\one_C\mid \Fc_1]>0$ on $A$.  Then we can construct a decreasing sequence of sets $(B_n)_{n\ge 0}$, $B_n\subset C$, such that  $0<\Er[\one_{B_n}\mid\Fc_1]\le 2^{-n}$ on $A$.
\end{lemma}
{\bf Proof}  The statement is obviously true for $n=0$ since we can take $B_0=C$.  We now proceed by induction and suppose the statement holds for $n$. So the set $B_n\subset A$ satisfies $0<\Er[\one_{B_n}\mid\Fc_1]\le 2^{-n}$ on $A$.  Clearly $A\subset \{\Er[\one_{B_n}\mid\Fc_1]>0\}$.  By assumption there is a set $D\subset B_n$ such that on $A\subset \{ \Er[\one_A\mid\Fc_1]>0\}$ we have
$$
0<\Er[\one_D\mid\Fc_1]<\Er[\one_{B_n}\mid\Fc_1].
$$
We now take
$$
B_{n+1}=\(D\cap \left\{\Er[\one_D\mid\Fc_1]\le \frac{1}{2}\Er[\one_{B_n}\mid\Fc_1]\right\}\)\cup
\((B_n\setminus D)\cap\left\{\Er[\one_D\mid\Fc_1]> \frac{1}{2}\Er[\one_{B_n}\mid\Fc_1]\right\}\).
$$
The set $B_{n+1}$ satisfies the requirements.
\begin{lemma} Let $C\in\Fc_2$ and let $h\colon\Omega\rightarrow [0,1]$ be $\Fc_1$ measurable.  Then there is a set $B\subset C$ such that $\Er[\one_B\mid\Fc_1]=h\,\Er[\one_C\mid\Fc_1]$.
\end{lemma}
{\bf Proof}  Let $B_0=\emptyset$.  Inductively we define for $n\ge 1$, classes $\Bc_n$ and sets $B_n\in \Bc_n$. For $n\ge 1$ let
$$
\Bc_n=\{  B_{n-1}\subset B\subset  C\mid B\in \Fc_2,\,\Er[\one_B\mid\Fc_1]\le  h\,\Er[\one_C\mid\Fc_1]\}.
$$
Let $\beta_n=\sup\{ \Pr[B]\mid B\in \Bc_n\}$ and take $B_n\in\Bc_n$ such that $\Pr[B_n]\ge (1-2^{-n})\beta_n$.  Clearly $B_n$ is non-decreasing and let  $B_\infty=\cup_n B_n$.  Obviously
$$
\Pr[B_\infty]\ge \limsup \beta_n\ge\liminf\beta_n\ge\lim \Pr[B_n]=\Pr[B_\infty].
$$
We claim that $\Er[\one_{B_\infty}\mid \Fc_1]=h\,\Er[\one_C\mid\Fc_1]$.  Obviously we already have that $\Er[\one_{B_\infty}\mid \Fc_1]\le h\,\Er[\one_C\mid\Fc_1]$.  If 
$\Pr\[ \Er[\one_{B_\infty}\mid \Fc_1] <  h\,\Er[\one_C\mid\Fc_1]  \]>0$ then $\Pr[B]<\Pr[C]$ and there must be $m\ge 1$ such that $\Pr\[ \Er[\one_{B_\infty}\mid \Fc_1] <  h\,\Er[\one_C\mid\Fc_1] -2^{-m}\]>0$.   The previous lemma allows to find $D\subset C\setminus B_\infty$, $\Pr[D]=\eta>0$,  such that $0<\Er[\one_D\mid\Fc_1]\le 2^{-m}$ on the set $\{\Er[\one_B\mid\Fc_1] <  h\,\Er[\one_C\mid\Fc_1] - 2^{-m}\}$ and zero elsewhere. The set $D\cup B_\infty$ is in all classes $\Bc_n$ and for $n$ big enough:
$$
\beta_n\ge\Pr[D\cup B_\infty] \ge \Pr[B_n]+\eta\ge (1-2^{-n})\beta_n +\eta\ge \beta_n+\eta-2^{-n}>\beta_n,
$$
yielding a contradiction.  So we must have $\Er[\one_{B_\infty}\mid \Fc_1]=h\,\Er[\one_C\mid\Fc_1]$.

\begin{rem}  The lemma above is a variant of Sierpi\'nski's theorem, \cite{Sier}. This theorem states that in an atomless probability space $(\Omega,\Ec,\Pr)$, for every set $A\in \Ec$ and every $0<t<1$, there is a set $B\subset A$ with $\Pr[B]=t\Pr[A]$. The usual proof --- presented in many probability courses --- uses the Axiom of Choice (AC).  A referee of \cite{del-com} pointed out that for many people  AC -- or Zorn's lemma -- is an extra assumption.  To prove Sierpi\'nski's theorem we only need the Axiom of Countable Dependent Choice, which is a countable form of the axiom of choice.  In analysis this is the axiom that is usually needed and used.   The proof above follows the approach given by Lorenc and Witu{l}a, \cite{LW}.
\end{rem}
\begin{lemma} There is an increasing family of sets $(B_t)_{t\in [0,1]}$ such that $\Er[\one_{B_t}\mid\Fc_1]=t$.  The sigma algebra $\Bc$, generated by the family $(B_t)_t$ is independent of $\Fc_1$.  The system $(B_t)_t$ can also be described as $B_t=\{U\le t\}$ where $U$ is a random variable that is independent of $\Fc_1$ and uniformly distributed on $[0,1]$.
\end{lemma}
{\bf Proof} The proof is a repeated use of the previous lemma where we take $h=1/2$. We start with $B_0=\emptyset, B_1=\Omega$.  Suppose that for the diadic numbers $k 2^{-n}, k=0,\ldots 2^n$ the sets are already defined.  Then we consider the set $B_{(k+1)2^{-n}}\setminus B_{k2^{-n}}$ and apply the previous lemma with $h=1/2$.  We get a set $D\subset B_{(k+1)2^{-n}}\setminus B_{k2^{-n}}$ with $\Er[\one_D\mid\Fc_1]=2^{-(n+1)}$.  We then define $B_{(2k+1)2^{-(n+1)}}=B_{k2^{-n}}\cup D$.  For non-diadic numbers $t$ we find a sequence of diadic numbers $d_n$ such that $d_n\uparrow t$.  Then we define $B_t=\cup_n B_{d_n}$.  This completes the construction.  Since the system $(B_t)_t$ is trivially stable for intersection, the relation $\Er[\one_{B_t}\mid\Fc_1]=t$ shows that the sigma algebra $\Bc$ generated by $(B_t)_t$, is independent of $\Fc_1$.  The construction of $U$ is standard.  At level $n$ we put $U_n=\sum_{k=1,\ldots 2^n} k2^{-n}\one_{B_{k2^{-n} }\setminus B_{(k-1)2^{-n}}}$.  $U_n$ then decreases to a random variable $U$ that satisfies the needed properties.
\begin{rem} After the  first version was made available, I got the remark that the paper \cite{SSWW} of Shen, J., Shen, Y.,  Wang, B., and  Wang, R.  contains similar concepts and results.\footnote{I thank Ruodu Wang for pointing out these relations and for the subsequent discussions we had on the topic.} In their notation they work with a measurable space $(\Omega,\Ac)$ on which they have a finite number of probability measures $\Qr_1,\ldots,\Qr_n$.\footnote{Their paper also considers an infinite number of measures but to clarify the relation between their paper and my approach, I only consider a finite number of measures.} They introduce
\begin{defi} The set $(\Qr_1,\ldots,\Qr_n)$ is conditionally atomless if there exists a dominating measure $\Qr$ (i..e $\Qr_k\ll \Qr$ for each $k\le n$) as well as a continuously distributed random variable $X$ (for the measure $\Qr$) such that the vector of Radon-Nikodym derivatives $\(\frac{d\Qr_k}{d\Qr}\)_{k}$ is independent of $X$.
\end{defi}
They then prove the following
\begin{prop} Are equivalent
\begin{enumerate}
\item $(\Qr_1,\ldots,\Qr_n)$ is conditionally atomless
\item in the definition we can take $\Qr=\frac{1}{n}(\Qr_1+\ldots+\Qr_n)$
\item $X$ can be taken as uniformly distributed over $[0,1]$.
\end{enumerate}
\end{prop}
There are several differences with my approach.  There is the technical difference that they suppose the existence of a continuously distributed random variable $X$.  In doing so they avoid the technical points between the more conceptual definition using conditional expectations and the construction of a suitable sigma-algebra with a uniformly distributed random variable.   A further difference is that they use a dominating measure that later can be taken as the mean of $(\Qr_1,\ldots,\Qr_n)$.  Of course their result together with the results here show that the definition of $(\Qr_1,\ldots,\Qr_n)$ being conditionally atomless, is equivalent to the statement that for the measure $\Qr_0=\frac{1}{n}(\Qr_1+\ldots+\Qr_n)$, the sigma algebra $\Ac$ is conditionally atomless with respect to the sigma-algebra generated by the Radon-Nikodym derivatives $\(\frac{d\Qr_k}{d\Qr_0}\)_k$.  In \cite{SSWW} it is also shown that one can take any strictly positive convex combination of the measures $(\Qr_1,\ldots,\Qr_n)$. Below we will show that this sigma-algebra in some sense has a minimal property, a result that clarifies the relation between the two approaches.  Before doing so, let us recall two easy results from introductory probability theory.
\begin{result}  For a given probability space $(\Omega,\Ac,\Qr)$ let us denote $\Nc=\{N\in\Ac\mid\Qr[N]=0\}$.  Suppose that a sub sigma-algebra $\Fc\subset\Ac$ is given and that $\Gc$, $\Fc\subset\Gc$, is another sub sigma-algebra which is included in the sigma-algebra generated by $\Fc$ and $\Nc$. Then for each $\xi\in L^1(\Omega,\Ac,\Qr)$
$$
\Er_\Qr[\xi\mid \Fc]=\Er_\Qr[\xi\mid \Gc] \quad\text{\as}
$$
\end{result}
\begin{result} With the notation in the previous exercise let $F\colon\Omega\rightarrow \Rr^n$ and $F'\colon\Omega\rightarrow \Rr^n$ be two vectors that are equal \as.  Let $\Fc$ be generated by $F$ and $\Gc$ be generated by $F'$.  Then $\Fc$ and $\Gc$ are equal up to sets in $\Nc$.  More precisely $\Gc$ is included in the sigma-algebra generated by $\Fc$ and $\Nc$ (and of course conversely), i.e. $\sigma(\Fc,\Nc)=\sigma(\Gc,\Nc)$.
\end{result}
\begin{prop} Let $\Qr_1,\dots,\Qr_n$ be probability measures on a measurable space $(\Omega,\Ac)$.  Let $\Qr_0$ denote  a convex combination of these measures $\Qr_0=\sum_k \lambda_k \Qr_k$ where each $\lambda_k > 0$. Let $f_k$ denote an $\Ac$ measurable version $\frac{d\Qr_k}{\Qr_0}$.  Let $\Qr$ be another dominating measure with $g_k$ an $\Ac$ measurable version of $\frac{d\Qr_k}{d\Qr}$.  Let $\Nc=\{N\in \Ac\mid \Qr_0[N]=0\}$.  Let $\Fc$ be generated by $f_k,k=1\dots n$ and let $\Gc$ be generated by $g_k,k=1\ldots n$.  Then $\Fc\subset \sigma(\Gc,\Nc)$
\end{prop}
{\bf Proof }  Clearly $\Qr_0\ll \Qr$ so let $h=\frac{d\Qr_0}{d\Qr}$. It is now immediate that $g_k= f_k h$ $\Qr$ \as. To see this, observe that the values of $f_k$ on $\{h=0\}$ do not matter. The functions $g_k$ and $h$ are $\Gc$ measurable since $h$ can be taken as $h=\sum_k \lambda_k g_k$. Then we define $f_k'= \frac{g_k}{h}$ on $\{h>0\}$ and $f_k'=0$ on $\{h=0\}$. This choice shows that the $f_k'$ are $\Gc$ measurable. It is immediate that $f_k=f_k'$ $\Qr_0$ \as.   The result now follows.

From the theorem it follows that the sigma-algebra augmented with the class $\Nc$ is the same for all strictly positive convex combinations.  The theorem shows that in the definition of conditionally atomless with respect to $\Fc$, we can also add the null sets $\Nc$ to $\Fc$.  To check that $\Ac$ is conditionally atomless with respect to a sigma-algebra $\Fc$ it is clear that the smaller $\Fc$, the easier it is to satisfy the condition. In my opinion the above clarifies the relation between this paper and \cite{SSWW}.
\end{rem}
\section{An Equivalent Definition}

As already mentioned in the previous section, the definition of being conditionally atomless is related to a similar statement for the kernel $K$.  
We suppose that the conditional expectation with respect to $\Fc_1$ is given by the kernel $K$.  
We have the following 
\begin{theo} If $\Fc_2$ is conditionally atomless with respect to $\Fc_1$ then for almost every $\omega\in \Omega$ the probability measure $K(\omega,.)$ is atomless on $\Fc_2$.  In case the sigma algebra $\Fc_2$ is generated by a countable family of sets, the converse holds, i.e. if for almost every $\omega\in\Omega$, the probability $K(\omega,.)$ is atomless on $\Fc_2$, then $\Fc_2$ is conditionally atomless with respect to $\Fc_1$.
\end{theo}
{\bf Proof } We first suppose that $\Fc_2$ is atomless with respect to $\Fc_1$.  According to the previous section there is an atomless sub sigma algebra $\Bc\subset \Fc_2$ that is independent of $\Fc_1$.  There is also a random variable $U$ which is independent of $\Fc_1$ and is uniformly distributed on $[0,1]$. Let $\Cc[0,1]$ be the space of real valued continuous functions on $[0,1]$, equipped with the sup norm.  This space is separable and so we can take a (sup-norm) dense sequence $(g_n)_{n\ge 1}$ in $\Cc[0,1]$.  For each $n\ge 1$ we have \as:
$$
\Er[g_n(U)\mid\Fc_1](\omega)=\Er[g_n(U)]=\int_0^1g_n(t)\,dt.
$$
So we have \as , say on $\Omega_n,\Pr[\Omega_n]=1$;
$$
\int K(\omega,d\tau)g_n(U(\tau))=\int_0^1g_n(t)\,dt.
$$
For  $\omega\in \cap_{n\ge 1}\Omega_n$ we have by density of the sequence $(g_n)_n$, for all $g\in \Cc[0,1]$:
$$
\int K(\omega,d\tau)g(U(\tau))=\int_0^1g(t)\,dt.
$$
This proves that \as the random variable $U$ is for $K(\omega,.)$ uniformly $[0,1]$ distributed.  That can only happen when $K(\omega,.)$ is atomless on $\Fc_2$.

We now prove the converse.  Suppose that $\Fc_2$ is not conditionally atomless with respect to $\Fc_1$.  In this case there is a set $A$ with $\Pr[A]>0$ such that for all $B\subset A$:
$$
\Pr\[0 < \Er[\one_B\mid \Fc_1] < \Er[\one_A\mid \Fc_1] \] =0.
$$
In order words, if $B\subset A$ then \as either $\Er[\one_B\mid \Fc_1]=0$ or $\Er[\one_B\mid \Fc_1]=\Er[\one_A\mid \Fc_1]$.  By definition of the kernel $K$, this means $K(\omega,B)=0$ or $K(\omega,B)=K(\omega,A)$ \as.  In other words for $B\subset A$:,
$$
\Pr\[\left\{ \omega\mid K(\omega,B)=0 \text{ or } K(\omega,B)=K(\omega,A)   \right\}     \] = 1.
$$
Since $\Fc_2$ is countably generated there is a countable Boolean algebra $\Ac\subset \Fc_2$ that generates $\Fc_2$. For each set $B\in \Ac$ we have that
$$
\Omega_B=\left\{ \omega\mid K(\omega,B\cap A)^2=K(\omega,A)\,K(\omega,B\cap A)  \right\},
$$
has measure $1$.  The set $\Omega'=\cap_{B\in\Ac}\Omega_B$ still has probability $1$. We claim that for each $\omega\in \Omega'$ and each $B\in \Fc_2$ we have that either $K(\omega, B\cap A)=0$ or $=K(\omega,A)$.  This means that for each $\omega\in \Omega'$ with $K(\omega,A)>0$, the measure $K(\omega,.)$ has $A$ as an atom, a contradiction to the hypothesis.  To show the claim we use a monotone class argument.  Let
$$
\Mc= \left\{ B\in \Fc_2\mid \text{ for each }\omega\in\Omega': K(\omega,B\cap A)^2=K(\omega,A)\,K(\omega,B\cap A)  \right\}.
$$
Clearly $\Ac\subset \Mc$ and it is obvious that $\Mc$ is a monotone class, meaning that it is stable for increasing countable unions and for decreasing countable intersections. It is well known that this implies $\Mc=\Fc_2$, completing the proof of the theorem.

\section{A Counterexample}
We now give a counterexample when $\Fc_2$ is not countably generated. The basic ingredient is the interval $[0,1]$ with its Borel sigma algebra $\Bc$ and the Lebesgue measure $m$.  We define $S=[0,1]$ and $\Omega=[0,1]\times (S\times [0,1])$.  The sigma algebra $\Fc_1$ is generated by the first coordinate and the Borel sigma algebra, $\Bc$, on $[0,1]$.  On $S\times [0,1]$ we put the sigma algebra defined as follows:
$$
\Ac=\left\{ B\mid \text{ there is a countable set } D \text{ and for } s\in D: B_s\in \Bc \text{ otherwise } B_s=\emptyset \text{ or }  B_s=[0,1]  \right\}.
$$
The sigma algebra $\Fc_2$ is the product sigma algebra $\Bc\otimes \Ac$.
For each $x\in [0,1]$ we define the kernel $K(x,C)$ as follows.  We first define the transition probability $k(x,B)$ for $B\in \Ac$:
$$
k(x,B)=\sum_{s\in S}\one_{\{x\}}(s)\,m(B_s).
$$
Then we define the kernel (defined on $\Omega$) as $K(\omega,.)=\delta_x\otimes k(x,.)$, where $\omega=(x,s,y)$ and where $\delta_x$ is the Dirac measure concentrated on the point $x$. The probability measure on $\Omega$ is constructed with $m$ and the transition kernel $k$:
$$
E\in \Bc, B\in \Ac\quad\Pr[E\times B]=\int_E m(dx) \, m(B_x)= m\(E\cap \{x\mid B_x=[0,1]\}\).
$$
For each $x \in [0,1]$ the kernel $K$ is atomless.\\
For $A\in \Fc_2$ we find putting $A_x=\{(s,z)\mid (x,s,z)\in A\}$ and $A_{x,x}=\{y\mid (x,x,y)\in A\}$:
$$
\Er[\one_A\mid \Fc_1](x) = k(x,A_x)=m(A_{x,x}),
$$
which is almost surely $0$ or $1$.  This makes it impossible that $\Fc_2$ is conditionally atomless with respect to $\Fc_1$.

\end{document}